\documentclass{amsart} 

\usepackage{amsmath,amsthm}     
\usepackage{graphicx}     
\usepackage{amsfonts} 

\usepackage{multicol}

\begin{document}

\begin{abstract}
We introduce \textsc{Partiti}, the puzzle that will run in the \textsc{Magazine} this year, and use the opportunity to recall some basic properties of integer partitions.
\end{abstract}

\author{Andr\'es Eduardo Caicedo}
\address{           
Mathematical Reviews\\    
416 Fourth Street\\ 
Ann Arbor MI, 48104
}
\email{aec@ams.org}
\thanks{\textbf{Andr\'es Eduardo Caicedo} earned his BSc in Mathematics from Universidad de los Andes, in Bogot\'a, Colombia, in 1996, and his PhD from the University of 
California, Berkeley, in 2003. His research centers on set theory. From 2003--2005, he was a research assistant at the Kurt G\"odel Research Center for Mathematical Logic of the 
University of Vienna. From 2005--2008 he was the Harry Bateman Research Instructor at the California Institute of Technology. He then joined the department of mathematics at Boise 
State University, where he remained until moving to Ann Arbor in 2015 as an associate editor at Mathematical Reviews. Andr\'es is married and has two kids. He enjoys comic books, 
Sudoku and Candy Crush.}
\urladdr{http://www-personal.umich.edu/~caicedo/}

\author{Brittany Shelton}
\address{
Albright College \\
Reading, PA 19612
}
\email{bshelton@albright.edu}
\thanks{\textbf{Brittany Shelton} earned her BSc from Montclair State University, in 2007, and her PhD from Lehigh University, in 2013. She is an assistant professor of 
mathematics at Albright College, Reading PA. Her research interests include algebraic and enumerative combinatorics. She never gives up an opportunity to combine two of her 
passions: mathematics and puzzles.}



\title{Of puzzles and partitions: Introducing Partiti}

\maketitle

In each of the five issues for 2017, readers of this \textsc{Magazine} found a \textsc{Pinemi} puzzle. Pinemi is the creation of Vietnamese puzzle enthusiast Thinh Van Duc Lai, who has 
also designed \textsc{Partiti}, the puzzle that will run through this year's issues. 

\section{Partiti}

Partiti is played on a $6\times 6$ grid in which each cell contains a positive integer. To play, place one or more digits into each cell in such a way that the digits in a cell sum to the 
indicated positive integer and no digit appears more than once in a cell or between cells that are adjacent or share a corner. 

The objective of the game can be described as finding unordered integer partitions of the given numbers consisting of distinct parts from $1,2,\dots,9$ (subject to the additional restriction 
that the partitions for contiguous cells should use different parts). Such an integer partition of $n$ consists of an increasing sequence of positive integers that sum to $n$. See Figure 
1 for an example.

\begin{figure}[h]
\centering
\includegraphics{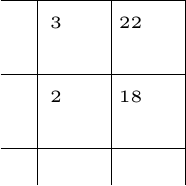}
\caption{The top-right corner of this month's puzzle. We can solve some of it by noting that the only partition of $2$ into distinct parts is $2$ itself, and the only such partitions of $3$ are 
$1+2$ and $3$, but the former is excluded, since we cannot use 2 again. Though more information is needed to see what numbers go into the other two cells, the reader may want to 
note that $22+18=1+4+5+6+7+8+9$, the sum of the remaining digits, so all available numbers should be used between these two cells.}
\end{figure}

The puzzle for this month is at the end of this note. In what follows, we present some basic properties of integer partitions, take a very brief detour through partitions of infinite sets, and 
conclude with a few words about Partiti's creator Thinh Lai.

\section{Integer partitions}

The study of partitions began with Euler. The number of integer partitions of $n$ is often denoted $p(n)$. Hardy and Ramanujan worked out an analytic formula for $p(n)$; the 
formula takes the form of an infinite series, and even just a few terms  produce remarkably accurate approximations.  As $n$ increases, $p(n)$ grows faster than a polynomial, but 
slower than any exponential $a^n$, $a>1$. More specifically, $p(n)\sim \frac1{4n\sqrt3}{e^{\pi\sqrt{2n/3}}}$, meaning that the quotient of these two expressions approaches one as 
$n$ approaches infinity. However, puzzlers need not worry, since the possible positive integers in the cells of Partiti are quite modest, the largest potential entry in a cell being 
$9+8+7+6+5+4=39$, which could only occur in a corner surrounded by $1$, $2$, and $3$.

More relevant than $p$ in this context is the number of partitions of $n$ into distinct parts, usually denoted $q(n)$. For example, $5=1+4=2+3$ are all the partitions of 
$5$ into distinct parts, and therefore $q(5)=3$. By convention, $q(0)=1$. The function $q$ has a somewhat more modest rate of growth than that of $p$, namely, $q(n)\sim 
\frac1{4\root4\of{3n^3}}e^{\pi\sqrt{n/3}}.$ The sequence $q(0),q(1),\dots$ is sequence A000009 in the OEIS \cite{OEIS}.

The first nontrivial result on this function $q$ is Euler's theorem from 1748 \cite{Euler} giving us that $q(n)$ coincides with the number of partitions of $n$ into (not necessarily distinct) 
odd parts. For example, $5=1+1+3=1+1+1+1+1$ are all such partitions of $5$ and, as predicted by the theorem, there are precisely $3$ of them. See Figure 2 for the beginning 
of Euler's work on integer partitions.

\begin{figure}[h] 
\centering
\includegraphics[height=45mm]{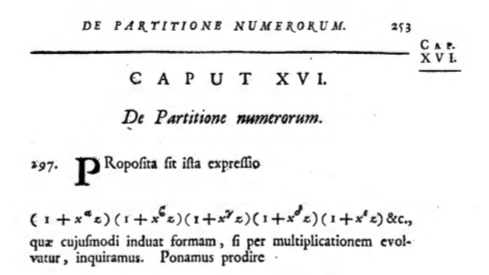}
\caption{The opening of chapter 16 of Euler's \textit{Introductio in Analysin Infinitorum} \cite{Euler}. The book lays the foundations of mathematical analysis. It also introduces the theory 
of integer partitions, in this chapter.}
\end{figure} 

There are several proofs of Euler's theorem. The one we proceed to sketch uses generating functions. It relies on observing that $\sum_{n=0}^\infty q(n)x^n$ can be represented as 
$\prod_{n=1}^\infty(1+x^n)$: At least formally, this product can be expanded by picking from each factor $1+x^n$ one of the two summands, with the understanding that in each 
product, all but finitely many times the summand $1$ is the chosen one. Grouping together like terms, the coefficient of $x^n$ in this expansion counts the number of ways the 
exponent $n$ can be formed as a sum of distinct terms, which is precisely $q(n)$. For example, note that the only products that result in an $x^5$ term are $1\cdot1\cdot1\cdot1\cdot 
x^5=x\cdot1\cdot 1\cdot x^4=1\cdot x^2\cdot x^3$, where in each product we have omitted the infinitely many remaining $1$s. 

Now, note that $\prod_{n=1}^\infty (1+x^n)=\prod_{n=1}^\infty\frac{1-x^{2n}}{1-x^n}=\prod_{n=1}^\infty\frac1{1-x^{2n-1}}$, the latter equality holding because all numerators cancel out 
and the only denominators that survive are the ones with odd degree. Expanding this product reveals that it is the generating function for partitions into odd parts:
{\small
$$\begin{aligned}
\prod_{n=1}^\infty\frac1{1-x^{2n-1}}&=\prod_{n=1}^\infty(1+x^{2n-1}+x^{2(2n-1)}+x^{3(2n-1)}+\cdots)\\
&=(1+x+x^{2\cdot 1}+\cdots)(1+x^3+x^{2\cdot 3}+\cdots)(1+x^5+x^{2\cdot 5}+\cdots)\cdots\\
&=(1+x+x^1 x^1+\cdots)(1+x^3+x^3x^3+\cdots)(1+x^5+x^5x^5+\cdots)\cdots\\
&=1+ x+ x^1x^1+(x^1x^1x^1+x^3)+(x^1x^1x^1x^1 +x^1x^3)+\cdots.
\end{aligned}$$}

The argument above can be readily formalized either in terms of formal power series expansions or in terms of ``genuine'' power series (upon arguing that the series involved 
converge for $|x|<1$).

Many other interesting results are known for $q$ and other partition functions, see \cite{AndrewsEriksson} for an introduction. These results are established by a wide variety of 
techniques, including combinatorial counting arguments, Ferrers diagrams, and others, and some are quite sophisticated, involving detailed analytical arguments, which entered the 
picture thanks to the joint work of Hardy and Ramanujan at the beginning of the twentieth century. 

\section{A small detour}

The first-named author cannot help but mention that some of his own work involves the study of partitions, in this case partitions of infinite sets. This is part of the area of set theory 
called the partition calculus. As a simple example of the sort of problems one considers here, readers familiar with the distinction between countable and uncountable sets may 
enjoy verifying the following: Suppose the set $\mathbb R$ of reals is partitioned into countably many pieces, $\mathbb R=\bigcup_{n=1}^\infty A_n$. Then at least one of the sets 
$A_n$ contains an infinite increasing sequence. Note the result fails for $\mathbb Q$ in place of $\mathbb R$ (we can split $\mathbb Q$ into countably many singletons). On the 
other hand, the result is not simply an artifact of $\mathbb R$ being uncountable (which ensures that one of the $A_n$ is also uncountable), since not every uncountable ordered set 
contains an infinite increasing sequence.

\section{About Partiti's creator}

We hope readers enjoy Partiti and the many other puzzles we anticipate seeing from Thinh. They can learn more about Thinh himself in a recent piece on his work by Will Shortz 
that ran in The New York Times \cite{Shortz}. Thinh's puzzle Bar Code appeared for 14 weeks in the Sunday Magazine section of The Times. 

Thinh has shown us a large variety of different puzzles of his own creation, all of a somewhat mathematical flavor. We asked him a few questions in preparation for this note. He 
shared with us that he solves all his puzzles manually, and uses the time it takes him to estimate their difficulty. He has surprised himself a few times with how hard some of his 
creations turned out to be. Although the name ``Partiti'' is probably self-explanatory, most of his puzzle names are inspired by Japanese puzzles, of which he confesses to be a big fan.
 
Thinh advises readers interested in creating their own mathematical puzzles that it is necessary to build a basic foundation and to read about numbers and logic puzzles. It took 
him five years to acquire this foundation himself. He indicates that some of the examples he has submitted to this \textsc{Magazine} are harder than the ones he has had featured in 
The New York Times, and hopes to publish a book of his own puzzles.

\section{Acknowledgment}
We would like to thank Thinh Lai for his enthusiasm and help with the preparation of this note.


\begin{thebibliography}{3}

\bibitem{AndrewsEriksson} G. Andrews and K. Eriksson, \textit{Integer partitions}. Cambridge University Press, Cambridge, 2004. 

\bibitem{Euler} L. Euler, \textit{Introductio in analysin infinitorum}. Apud Marcum-Michaelem Bousquet \& socios, Laussane, 1748. https://archive.org/details/%
bub$\underline{\ }$gb$\underline{\ }$jQ1bAAAAQAAJ.

\bibitem{OEIS} Sequence A000009. The on-Line encyclopedia of integer sequences. https://oeis.org/A000009.

\bibitem{Shortz} W. Shortz, Bar Code: a new infatuation poised for a puzzle craze, \textit{The New York Times}, June 12, 2017. https://nyti.ms/2sfDb1f.

\end{thebibliography}
\end{document}